\documentclass[preprint,12pt]{elsarticle}

\usepackage{amssymb}
\usepackage{amsmath}
\usepackage[ruled,vlined]{algorithm2e}
\usepackage{multirow}

\journal{Russian Journal of Numerical Analysis and Mathematical Modeling}

\begin{document}
	
	\begin{frontmatter}
		
		%% Title, authors and addresses
		
		%% use the tnoteref command within \title for footnotes;
		%% use the tnotetext command for theassociated footnote;
		%% use the fnref command within \author or \address for footnotes;
		%% use the fntext command for theassociated footnote;
		%% use the corref command within \author for corresponding author footnotes;
		%% use the cortext command for theassociated footnote;
		%% use the ead command for the email address,
		%% and the form \ead[url] for the home page:
		%% \title{Title\tnoteref{label1}}
		%% \tnotetext[label1]{}
		%% \author{Name\corref{cor1}\fnref{label2}}
		%% \ead{email address}
		%% \ead[url]{home page}
		%% \fntext[label2]{}
		%% \cortext[cor1]{}
		%% \address{Address\fnref{label3}}
		%% \fntext[label3]{}
		
		\title{Comparison of nonlinear solvers within continuation method for steady-state variably saturated groundwater flow modeling}
		
		%% use optional labels to link authors explicitly to addresses:
		%% \author[label1,label2]{}
		%% \address[label1]{}
		%% \address[label2]{}
		
		\author[label1,label2]{Denis Anuprienko}
		%\author[label1,label2]{Ivan Kapyrin}
		\address[label1]{Marchuk Institute of Numerical Mathematics, Russian Academy of Sciences}
		\address[label2]{Nuclear Safety Institute, Russian Academy of Sciences}
		
		%\address{Address}
		
		\begin{abstract}
			Nonlinearity continuation method, applied to boundary value problems for steady-state Richards equation, gradually approaches the solution through a series of intermediate problems.  Originally, the Newton method with simple line search algorithm was used to solve the intermediate problems. In this paper, other solvers such as Picard and mixed Picard-Newton methods are considered, combined with slightly modified line search approach. Numerical experiments are performed with advanced finite volume discretizations on model and real-life problems.
		\end{abstract}

		\begin{keyword} Groundwater flow \sep Unsaturated conditions \sep Vadose zone \sep Richards equation \sep Newton method \sep Line search \sep Continuation
			%% keywords here, in the form: keyword \sep keyword
			
			%% PACS codes here, in the form: \PACS code \sep code
			
			%% MSC codes here, in the form: \MSC code \sep code
			%% or \MSC[2008] code \sep code (2000 is the default)
			
		\end{keyword}
		
	\end{frontmatter}
	
	%% \linenumbers
	
	%% main text
	\section{Introduction}
	\label{}
	Mathematical modeling of groundwater flow is important in various problems in hydrogeology. Many of these problems consider near-surface zone, where soil pores are only partially saturated with water. Prediction of groundwater flow in such conditions requires solution of nonlinear partial differential equation, the Richards equation \cite{richards1931capillary, bear2010modeling}. Some analytical and semi-analytical solutions \cite{parlange1999analytical} are available for the Richards equation, but solution is usually found numerically. Difficulties encountered in process are related to nonlinear systems resulting from discretization of the Richards equation both in transient and steady-state cases with many approaches having been developed (see, for example, \cite{celia1990general, farthing2003efficient, diersch1999primary, pop2004mixed}). This paper focuses on steady-state case and, more precisely, on a special solution technique applied recently for these problems -- the nonlinearity continuation method \cite{anuprienko2021nonlinearity}. Unlike commonly used pseudo-transient method \cite{farthing2003efficient, fowler2005pseudo, anuprienko2018modeling}, in which the system evolves in time until reaching steady-state, the nonlinearity continuation works directly with the steady-state equation. As it has been reported \cite{farthing2003efficient, paniconi1994comparison}, straightforward application of standard nonlinear solvers such as Newton method may fail due to lack of good initial guess. The nonlinearity continuation method obtains initial guess through solution of a sequence of intermediate problems with increasing nonlinearity, starting with a simple linear problem. So far, Newton method with simplest line search was used for solution of the problems in the sequence \cite{anuprienko2021nonlinearity}. In this paper, other approaches such as Picard and mixed Picard-Newton solvers are tested along with a line search technique based on Armijo rule. It is noteworthy to mention that the nonlinearity continuation procedure can be considered in predictor-corrector terms \cite{allgower2003introduction}. In its current form, the procedure has trivial predictor and the corrector is application of a nonlinear solver. Therefore, the scope of this article is comparison of different correctors. Sophisticated predictors are to be examined further.
	\par The paper is organized as follows. The second section contains mathematical description of steady-state groundwater flow in variably saturated conditions. The third section describes numerical solution approach and nonlinearity continuation method for arising nonlinear systems. The fourth section focuses on various nonlinear solvers within the nonlinearity continuation method. The fifth section contains results of numerical experiments. Conclusions are provided in the end.

	\section{Governing equations}
	\label{}
	Steady-state groundwater flow in variably saturated conditions is governed by steady-state version of Richards equation \cite{richards1931capillary, bear2010modeling}:
	
	\begin{equation}\label{eq:Richards}
 -\nabla \cdot \left(K_r(\theta)\mathbb{K}\nabla \left(\psi + z\right)\right) = Q.
	\end{equation}
	Here the following variables are used:
	\begin{itemize}
		\item $\psi$ -- capillary pressure head;
		\item $\theta(\psi)$ -- volumetric water content in medium;
		\item $\mathbb{K}$ -- hydraulic conductivity tensor, a 3$\times$3 s.p.d. matrix;
		\item $K_r(\theta)$ -- relative permeability for water in medium;
		\item $Q$ -- specific sink and source terms.
	\end{itemize}
	
	In this paper, the steady-state Richards equation \eqref{eq:Richards} is considered in terms of hydraulic head $h = \psi + z$ (here $z$ is the vertical coordinate). Since $K_r(\theta) = K_r(\theta(\psi)) = K_r(\theta(h - z))$, relative permeability can be represented as a function of hydraulic head $h$ and for simplicity denoted as $K_r(h)$. Then, one can formulate steady-state Richards equation in the following form:
	\begin{equation}\label{eq:RichStat}
	-\nabla \cdot \left(K_r(h)\mathbb{K}\nabla h\right) = Q.
	\end{equation}
	
	\par
	
	In addition to Richards equation, constitutive relationships between water pressure head $\psi$, water content $\theta$ and relative permeability $K_r$ are required. In this paper, two models are considered.
	
	\subsection{Van Genuchten -- Mualem model}
	This widely used model is based on nonlinear functions for $\theta(\psi)$ and $K_r(\theta)$ proposed by van Genuchten \cite{van1980closed} and Mualem \cite{mualem1976new}:
	\begin{equation}\label{eq:vgm_theta}
	\theta(\psi) = \theta_r + \frac{\theta_s - \theta_r}{(1 + |\alpha \psi|^ n)^ m},
	\end{equation}
	
	\begin{equation}\label{eq:vgm_Kr}
	K_r(\theta) = S_e^{1/2} \cdot \left(1 - \left(1 - S_e^{1/m}\right)^m\right)^2,
	\end{equation}
	where $\theta_s$ is the water content at full saturation, $\theta_r$ is the residual water content, $S_e = (\theta - \theta_r)/(\theta_s - \theta_r)$ is the effective saturation and $n > 1$, $m = 1 - 1/n$ and $\alpha$ are parameters of the model.

	\subsection{Unconfined flow model}
	A much simpler model, so-called unconfined flow model uses piecewise linear dependencies proposed in \cite{anuprienko2018modeling} and is designed specifically for numerical solution with the finite volume method. These functions depend on spatial discretization, possibly varying from cell to cell.
	Namely, the dependence of water content $\theta_E$ in a cell $E$ on hydraulic head $h_E$ in that cell is defined as follows
	\begin{equation}\label{eq:unconf_theta}
	\theta_E(h_E) = 
	\begin{cases}
	\phi, h_E > h_{E,\max},\\
	\phi \cdot \frac{h_E - h_{E,\min}}{h_{E,\max} - h_{E,\min}}, h_{E,r} < h_E \le h_{E,\max},\\
	\phi \cdot (\alpha_{\phi} - \alpha_{\theta}(h_{E,r} - h_E)), h_E \le h_{E,r},
	\end{cases}
	\end{equation}
	where $\phi$ is the medium porosity, $h_{E,\max}$ and $h_{E,\min}$ are the maximal and minimal vertical node coordinates of the cell $E$ (computed during numerical solution procedure)
	and $h_{E,r}$ is such that water content calculated by the second linear part in \eqref{eq:unconf_theta} is equal to $\phi \cdot \alpha_{\phi}$, namely,
	\begin{equation}
	h_{E,r} = h_{E,\min} + \alpha_{\phi}(h_{E,\max} - h_{E,\min}).
	\end{equation}
	\par
	Relative permeability for a cell $E$ is assumed to be equal to the saturation:
	\begin{equation}\label{eq:unconf_Kr}
	K_{r,E}(h_E) = S(h_E) = \frac{\theta_E(h_E)}{\phi}.
	\end{equation}
	\par 
	Compared to van Genuchten and Mualem dependencies \eqref{eq:vgm_theta} and \eqref{eq:vgm_Kr}, the constitutive relations \eqref{eq:vgm_theta} and \eqref{eq:vgm_Kr} demand only porosity and saturated hydraulic conductivity as media parameters. Model parameters $\alpha_s$ and $\alpha_\theta$ are assigned small values only to assure the nonnegativity of moisture content, and are the same for all media present in the domain.  
	
	%% The Appendices part is started with the command \appendix;
	%% appendix sections are then done as normal sections
	%% \appendix

	\section{Numerical solution with the nonlinearity continuation method}

	Discretization of the steady-state Richards equation \eqref{eq:RichStat} is performed with the cell-centered finite volume method. Unstructured grids used consist of either triangular prisms with occurrence of other polyhedra in some cases or of hexahedra with octree-based refinement and cut cells on boundaries \cite{plenkin2015adaptive}.
	Anisotropic and highly heterogeneous conductivity tensors $\mathbb{K}$ are common in real-life hydrogeological problems, and therefore advanced discretization techniques should be used. Finite volume schemes are defined by approximation of flux $-K_r(h)\mathbb{K}\nabla h$ across cell face. In this paper, the following options are considered:
	\begin{itemize}
		\item conventional linear two-point flux approximation (TPFA);
		\item linear multipoint flux approximation  (MPFA-O) \cite{aavatsmark1998discretization};
		\item nonlinear monotone two-point flux approximation (NTPFA-B) \cite{terekhov2017cell};
		\item nonlinear multipoint flux approximation with discrete maximum principle (NMPFA-B)  \cite{terekhov2017cell}.
	\end{itemize}
    The relative permeability $K_r(h)$ for a face can be either calculated as a half-sum of values from the two neighboring cells (central approximation) or taken from the cell with greater water head value (upwind approximation). In the author's experience, central approximation in general results in better convergence of nonlinear solvers since the structure of arising system does not change during iterations. Still, central approximation can lead to oscillations in saturation and in some cases to related convergence problems. In this paper, central approximation is used in all cases.
\par
	Discretization leads to a system of nonlinear equations which takes the form
	
	\begin{equation}\label{eq:NonlinSys}
	F(h) \equiv A(h)h-b(h) = 0,
	\end{equation}
	where $A(h)$ is a matrix with solution-dependent coefficients.
	\par
	The system is then solved by nonlinearity continuation method \cite{anuprienko2021nonlinearity}. Methods of this type \cite{allgower2003introduction}, when applied to numerical solution of nonlinear partial differential equations (PDEs) \cite{continCFD, continMechAdap, continMechNewt}, approach solution incrementally through a sequence of auxiliary problems. In case of steady-state variably saturated groundwater flow, nonlinearity of the governing PDE lies in the relative permeability function $K_r(h)$. This function is then replaced by a function depending on \textit{continuation parameter} $q$:
	\begin{equation}\label{eq:ContFun}
	\mathcal{K}(h,q):~~\mathcal{K}(h,0)\equiv 1,~~\mathcal{K}(h,1)\equiv K_r(h),
	\end{equation}
	of which two options are tested: \textit{linear} function $\mathcal{K}_{lin}(h,q) = 1 + q\cdot (K_r(h) - 1)$ and  \textit{power} function $\mathcal{K}_{pow}(h,q) = (K_r(h))^q$.
	\par
	The equation \eqref{eq:RichStat} becomes equation 
	
	\begin{equation}\label{ContEq}
	-\nabla \cdot \left(\mathcal{K}(h,q)\mathbb{K}\nabla h\right) = Q.
	\end{equation}
	\par
	From 
	%$\mathcal{K}(h,q)$ definition 
	\eqref{eq:ContFun} it is clear that with $q = 0$ the equation \eqref{ContEq} is a linear equation and with $q = 1$ it is the original steady-state Richards equation \eqref{eq:RichStat}. The solution of the original equation is then found in the following way. First, a linear problem with $q = 0$ is solved (it represents fully saturated case). Note that solution of the linear problem still requires application of a nonlinear solver in case of nonlinear finite volume schemes. Then the obtained solution is used as initial guess in nonlinear solver for problem with $q = q_1 > 0$. Then the new solution is again used as initial guess for problem with $q = q_2 > q_1$. These steps are repeated until the problem with $q = 1$ (the original one) is solved.

	\section{Nonlinear solvers within the continuation method}
	The nonlinear system \eqref{eq:NonlinSys} is solved in an iterative way. At iteration $k+1$ the update vector $\Delta h$ is found and then the update procedure is performed:
	\begin{equation}\label{eq:Update}
	h^{k+1} = h^k + \Delta h.
	\end{equation}
	\par Various solvers differ in the way of obtaining $\Delta h$. Iterations continue until one of the following conditions is satisfied:
	\begin{equation*}
	||F(h^{k+1})||_2 < \varepsilon_{rel} \cdot ||F(h^0)||_2
	\end{equation*}
	or
	\begin{equation*}
	||F(h^{k+1})||_\infty < \varepsilon_{abs}.
	\end{equation*}
	\par
	Iterations also stop if maximal allowed number $nit_{\max}$ is reached or divergence tolerance is exceeded, which is usually set to $\varepsilon_{div} = 10^{15}$.
	\subsection{Newton method}
	Newton method is widely used for solution of nonlinear systems. At iteration $k+1$, the update vector $\Delta h$ is obtained by solving the following linear system:
	\begin{equation}
		J(h^k)\Delta h = -F(h^k),
	\end{equation}
    where $J = [\partial F_i / \partial h_j]_{ij}$ is the Jacobian matrix.
    \par
    Quadratic convergence makes Newton method very attractive. However, some requirements are imposed on smoothness of the function $F(h)$ and initial guess. In applications for discretized nonlinear PDEs, lack of initial guess is one of the main reasons for convergence difficulties. The whole nonlinearity continuation procedure is a way to improve initial guess, but simpler ways also exist, including application of a simpler solver at first iterations.

    \subsection{Picard method}
    This method expoits special structure of function $F(h)$ given by equation \eqref{eq:NonlinSys}. At iteration $k+1$ the new head value $h^{k+1}$ can be obtained simply by solving the following linear system:
    \begin{equation}\label{eq:UsualPicard}
    	A(h^k) h^{k+1} = b(h^k).
    \end{equation}
    In order to allow unified approach to both Newton and Picard methods, by subtracting $A(h^k)h^k$ from both sides of \eqref{eq:UsualPicard} the Picard iteration is rewritten as  
    \begin{equation}
    	A(h^k)\Delta h = -F(h^k),
    \end{equation}
    followed by update procedure \eqref{eq:Update}.
    \par
    The Picard method converges only linearly and requires operator $A^{-1}(h)b(h)$ having contraction properties \cite{kelley1995iterative}. However, it does not need derivatives calculation and may be considerably cheaper on per-iteration basis compared to Newton (depending on implementation). The linear convergence is the main drawback for the Picard method as a standalone solver, but lesser sensitivity to initial guess makes it a good candidate to be used in a mixed solver as a way to improve initial guess.
    \par It should be also noted that finite volume schemes with nonlinear flux aproximations rely on application of Picard solver to prove their important properties such as solution monotonicity or discrete maximum principle satisfaction \cite{kapyrin2007family, vassilevski2020parallel}. In case of these schemes the Jacobian matrix $J$ may be more dense than the matrix $A$.

    \subsection{Mixed Picard-Newton solver}
    In this solver, a number $nit_{pic}$ of Picard iterations is done before switching to Newton. This approach is known to result in better convergence of the Newton method since the initial guess is improved by Picard iterations \cite{paniconi1994comparison}.
    
    \subsection{Line search technique with Armijo rule}  
    Sometimes the update $\Delta h$ in \eqref{eq:Update} may appear too large, and it is scaled by some parameter $0 < \omega \le 1$. This scaling is called \textit{relaxation} \cite{paniconi1994comparison}. \textit{Line search} is a procedure in which this parameter is chosen to ensure that $h^k + \omega\cdot\Delta h$ is a better approximation of solution than $h^k$.
    Previously, the simplest version of line search was considered \cite{anuprienko2021nonlinearity} which simply required $||F(h^k + \omega\cdot\Delta h)||_2 \le ||F(h^k)||_2$. It is known, however, that this approach may lead to oscillations without convergence \cite{kelley1995iterative} and such behavior was sometimes observed while working on \cite{anuprienko2021nonlinearity}. That is the reason to use more sophisticated line search based on Armijo rule \cite{armijo1966minimization, kelley1995iterative}, which already has been applied for Richards equation \cite{farthing2003efficient, jones2001newton}. In this approach, condition for acceptance of $\omega$ is strengthened and the resulting line search procedure at iteration is described in algorithm \ref{alg:Armijo}.
    
    	\begin{algorithm}[H]\label{alg:Armijo}
    	\SetAlgoLined
    	$\omega = 1$\;
    	\While{$\omega > \omega_{min}$}{
    		\eIf{$||F(h^k + \omega \Delta h)||_2 < (1-\alpha\omega)\cdot||F(h^k)||_2$}{
    			$h^{k+1} = h^k + \omega \Delta h$\;
    			%                    $||r_{k+1}||_2 = ||r||_2$\;
    			%                  $||r_{k+1}||_{\infty} = ||r||_{\infty}$\;
    			break\;
    		}{
    			$\omega = \gamma \cdot \omega$\;
    		}
    	}
    	\If{$\omega < \omega_{min}$}{
    		line search failed\;
    	}
    	\caption{Line search based on Armijo rule}
    \end{algorithm}
    Here $\alpha$ is a parameter determining sufficient decrease of $||F(h)||_2$ (taken $\alpha = 10^{-4}$ as in \cite{kelley1995iterative})), $\gamma = 0.25$ is the decreasing factor for $\omega$.
	
	\subsection{Strategies of line search application}\label{sec:LSS}
	A number of parameters regarding algorithm \ref{alg:Armijo} should be chosen. First, minimal value of $\omega$, $\omega_{min}$, should be specified. In \cite{anuprienko2021nonlinearity}, $\omega_{min}$ is chosen such that 7 iterations of $\omega$ refinement can be done, which leads to $\omega_{min} \approx 6\cdot10^{-5}$. During calculations for this paper, however, it was noted that sometimes smaller values  may be helpful, and 10 iterations of $\omega$ refinement were allowed. Still, very small values of $\omega_{min}$ may lead to a great number of iterations, and it may be cheaper overall to stop the solver and solve the problem in another few additional continuation steps.
	\par
	Another topic is the fact that nonlinear solvers often exhibit growth of $||F||$ at several first iterations and then still converge successfully. Application of line search at these first iterations may cause more harm than good. A solution of this problem proposed for Newton method in \cite{anuprienko2021nonlinearity} is not to apply line search for the first 5 iterations. This works well in cases when $||F||$ does not grow fast enough at these iterations. In many cases, especially for complex real-life problems, $||F||$ may exhibit sharp increase right from the beginning, exceeding divergence tolerance. On the other hand, in these same cases line search still fails at the first iterations. A solution may be to perform a number of initial iterations, $nit_{nls}$ ($nls$ = no line search), with fixed relatively small $\omega$ before switching to line search. These first iterations may coincide with initial Picard iterations if the mixed solver is used, and in this work it is done so.
	
	\section{Numerical experiments}
	This section includes a number of tests comparing approaches described above on model and real-life problems. 
	\subsection{Implementation}
	Numerical solution is implemented in GeRa software \cite{kapyrin2015integral, gera-site} which was originally designed for modeling of groundwater flow and contaminant transport for safety assessment of radioactive waste repositories and nowadays includes a wider spectrum of physical process models. GeRa is written in C++ and its computational core is based on INMOST \cite{inmost-site, vassilevski2020parallel}, an open-source platform which provides a set of tools for numerical modeling on unstructured meshes. One of the most useful features of INMOST used in this work is its automatic differentiation capability which is used to construct the Jacobian matrix in the Newton method and matrix $A$ in the Picard method. Internal INMOST linear solver \texttt{Inner\_MPTILUC} \cite{vassilevski2020parallel, terekhov2020parallel}, a version of BiCGStab with second order Crout-ILU \cite{konshin2015parallel} with inverse-based condition number estimation \cite{li2003crout} and maximum product transversal reordering \cite{duff1999design} as preconditioner, is used.
	\par Of the finite volume schemes used, TPFA and MPFA-O are internally implemented within GeRa, while NTPFA-B and NMPFA-B are taken from external package developed by Terekhov \cite{terekhov2017cell}. It should be noted that all schemes besides MPFA-O are implemented so that the unknowns are cell-centered only and therefore result in systems with size equal to number of mesh cells. MPFA-O includes unknowns on Neumann boundaries and therefore produces larger systems.
	
	\subsection{Problems description}
	The first problem is a model problem describing groundwater flow in a square dam (taken from \cite{polubarinova2015theory}). The dam problem is used as a verification test for the unconfined flow model and as a performance benchmark \cite{anuprienko2021nonlinearity, anuprienko2018modeling}. In this paper, that problem is slightly modified to exclude seepage boundary condition and make it more challenging for spatial discretization. The dam is represented by a 10 m $\times$ 10 m square domain. The domain is composed of a homogeneous material, this time with a full hydraulic conductivity tensor $\mathbb{K} = R(\pi/6)\cdot diag\{K_0,~10\cdot K_0\}\cdot R(-\pi/6)$, where $R$ represents rotation matrix and $K_0$ = 0.864 m/day. The left boundary has constant hydraulic head value $h = 10$ m, the right boundary has constant hydraulic head value $h = 2$ m up to the point where $z = 2$ m. Elsewhere boundaries are impermeable. The problem was considered in three-dimensional setting in the following way: only one layer of cells was present along third dimension, and all introduced boundaries were considered impermeable.
	\begin{figure}
		\centering
		\includegraphics[width=0.4\textwidth]{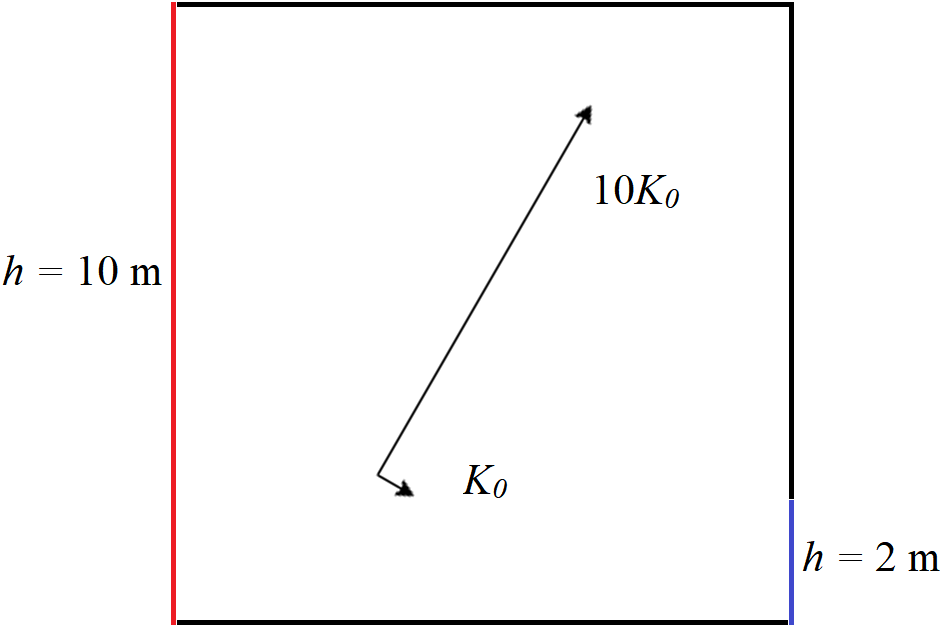}
		\caption{The dam domain, boundary conditions and hydraulic conductivity representation}\label{pic:dam_tensor}
	\end{figure}
	\par
	The second problem describes groundwater flow in a real site of approximately 64 km$^2$ featuring 7 different media forming 3 geological layers. All media have diagonal anisotropic hydraulic conductivity tensor of the form $\mathbb{K} = \text{diag}\{K, K, 0.1K\}$ with K varying from 0.011 to 4.76 m/day. This problem is used for performance benchmarking in GeRa and was also mentioned in \cite{anuprienko2018modeling} and \cite{anuprienko2021nonlinearity}. In this paper, the problem is modified to exclude lakes and rivers.
	% and replace non-uniform rainfall recharge distribution at the top boundary with the uniform one.
	 Originally, the problem uses unconfined flow model \eqref{eq:unconf_theta} -- \eqref{eq:unconf_Kr}, which turned out to be the reason for relatively good convergence of solvers \cite{anuprienko2021nonlinearity}, and its modification using the van Genuchten -- Mualem model is also considered.
	
	\subsection{Unconfined flow model: modified dam test}
	This first series of tests was performed on Cartesian cubic grids. Although structurally simple, these grids appear non-$\mathbb{K}$-orthogonal since $\mathbb{K}$ is rotated. Still, as the domain is geometrically trivial and homogeneous, it is fair to expect convergence of the nonlinearity continuation method  in 1 step. Therefore, solvers were tested in their ability to ensure such convergence. The tests were performed on Cartesian grids with 400 and 6400 cells (cell cize 0.5 and 0.125 m, respectively). The following solver parameters were used: $\varepsilon_{rel} = 10^{-5},~~\varepsilon_{abs} = 10^{-6},~~nit_{\max} = 50$. 
	\par
\begin{figure}
		\centering
		\includegraphics[width=0.4\textwidth]{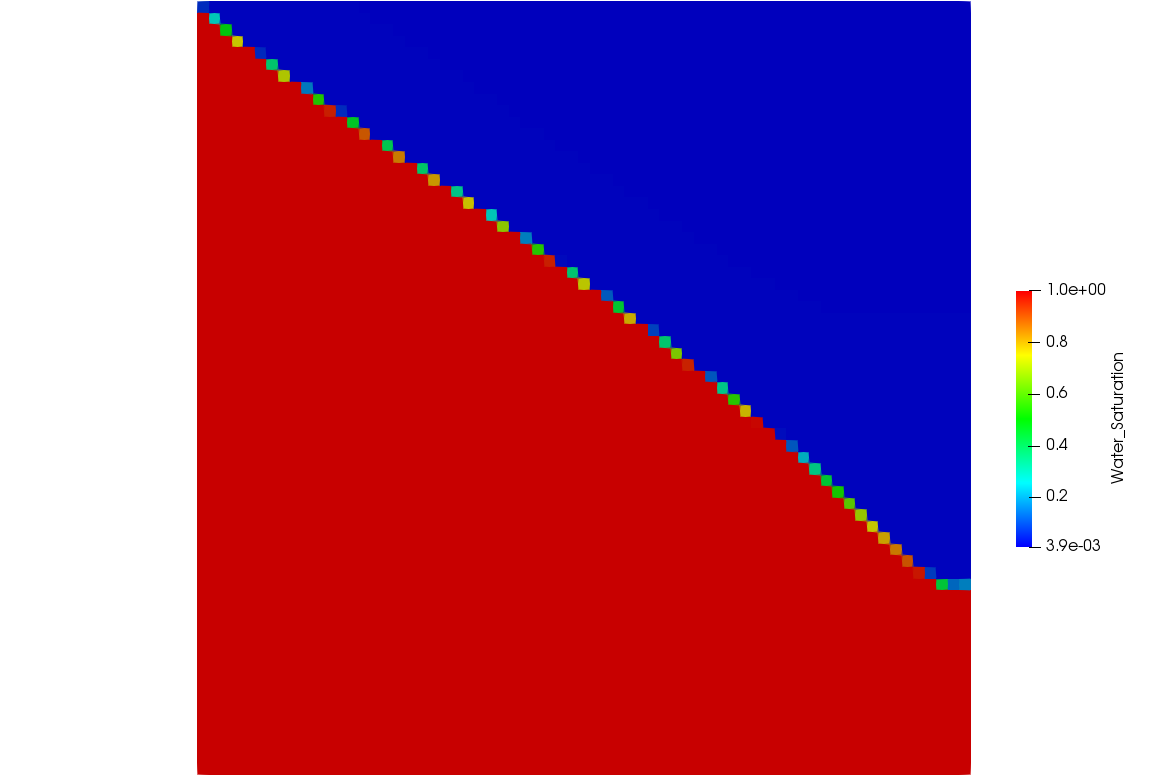}
		\includegraphics[width=0.4\textwidth]{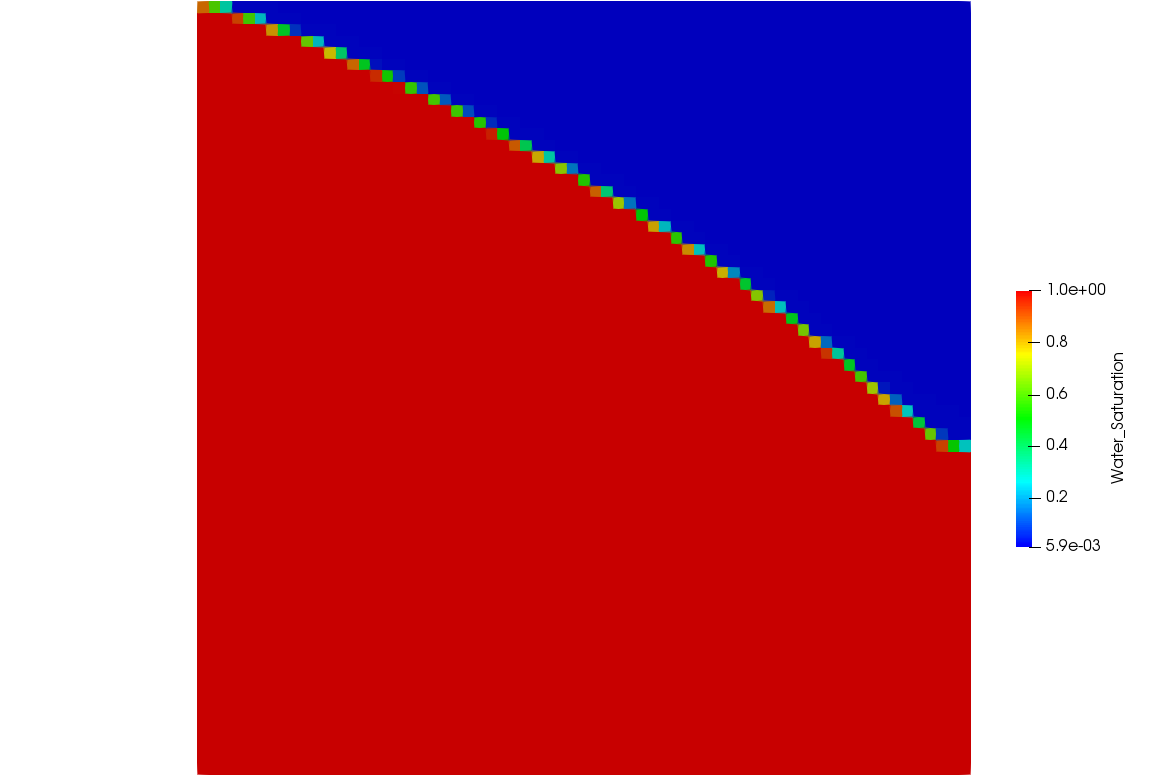}
		\caption{Saturation distributions for the modified dam problem produced by NMPFA-B and TPFA schemes, 6400 cells}\label{pic:dam_rot_sat}
	\end{figure}

\begin{figure}
		\centering
		\includegraphics[width=0.4\textwidth]{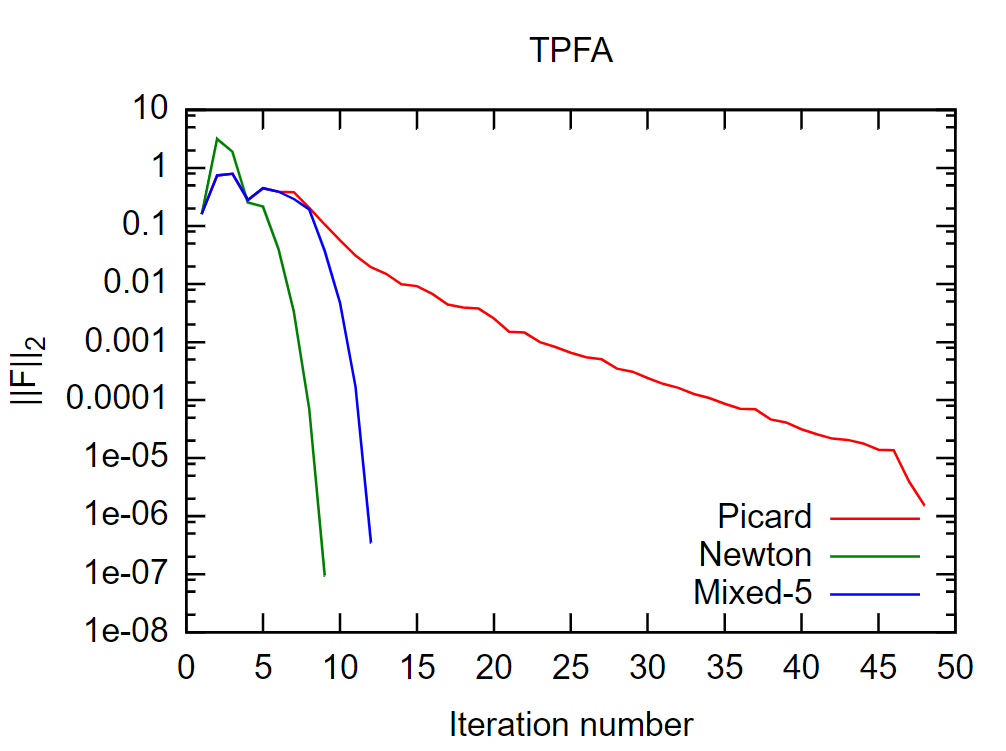}
		\includegraphics[width=0.4\textwidth]{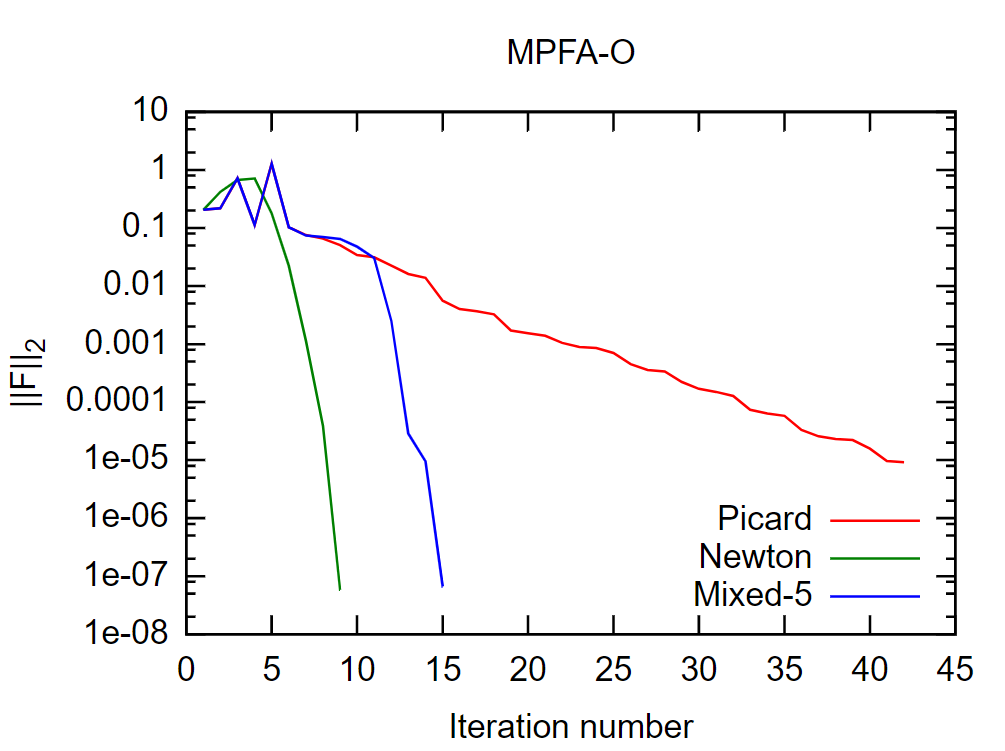}
		\includegraphics[width=0.4\textwidth]{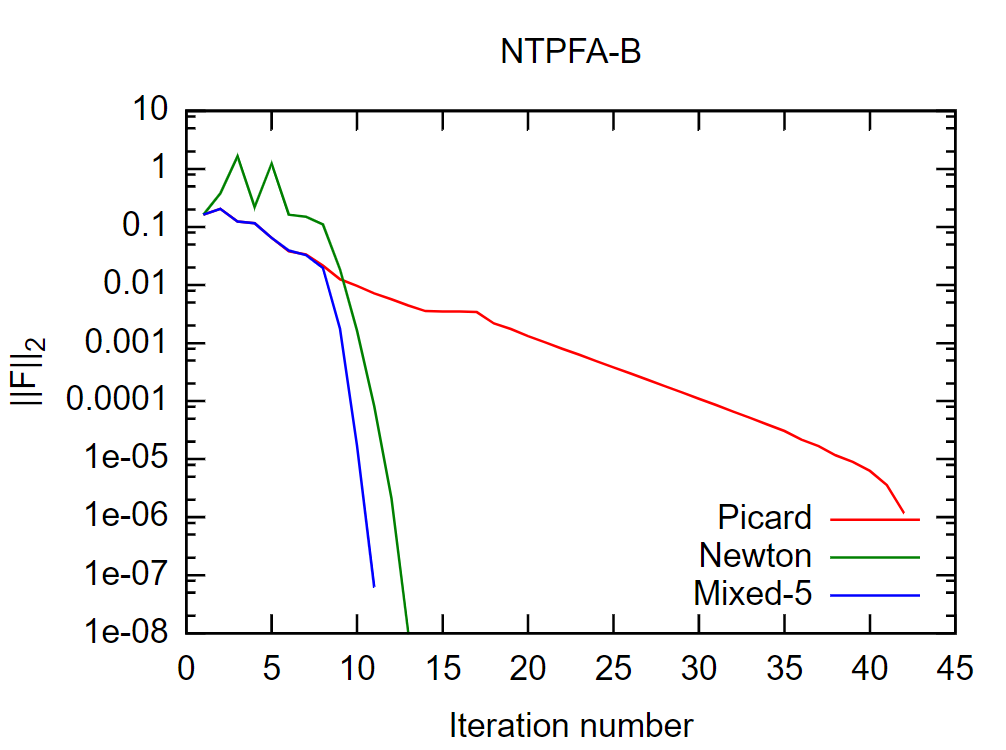}
		\includegraphics[width=0.4\textwidth]{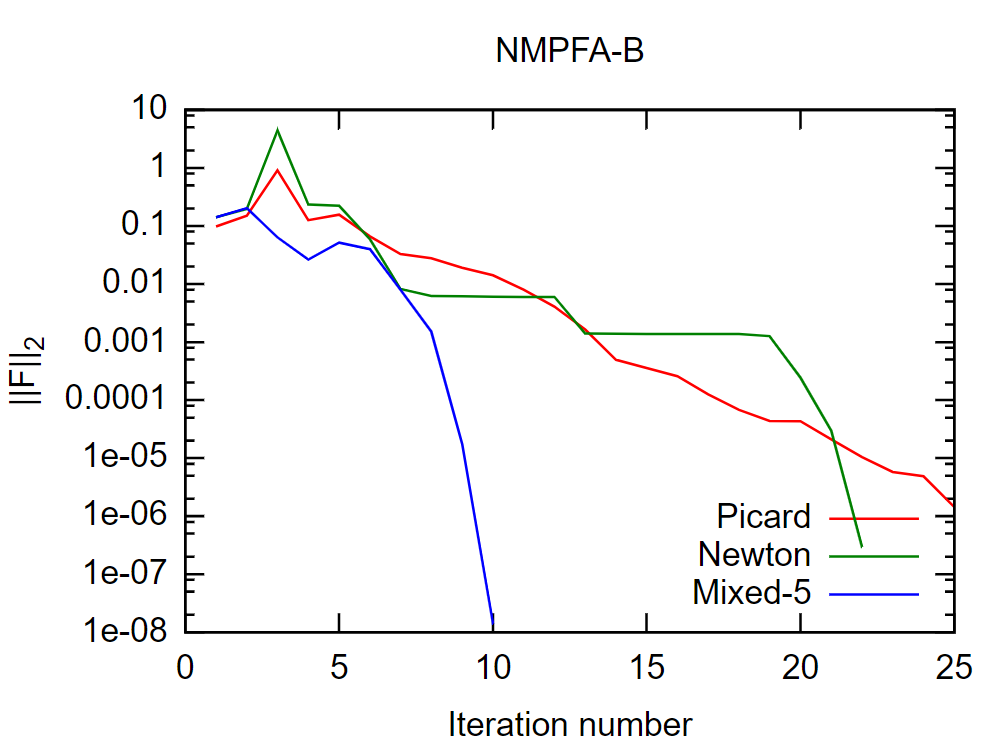}
		\caption{Convergence for the modified dam problem on Cartesian grid, 400 cells}\label{pic:dam_rot_400}
	\end{figure}

\begin{figure}
\centering
\includegraphics[width=0.4\textwidth]{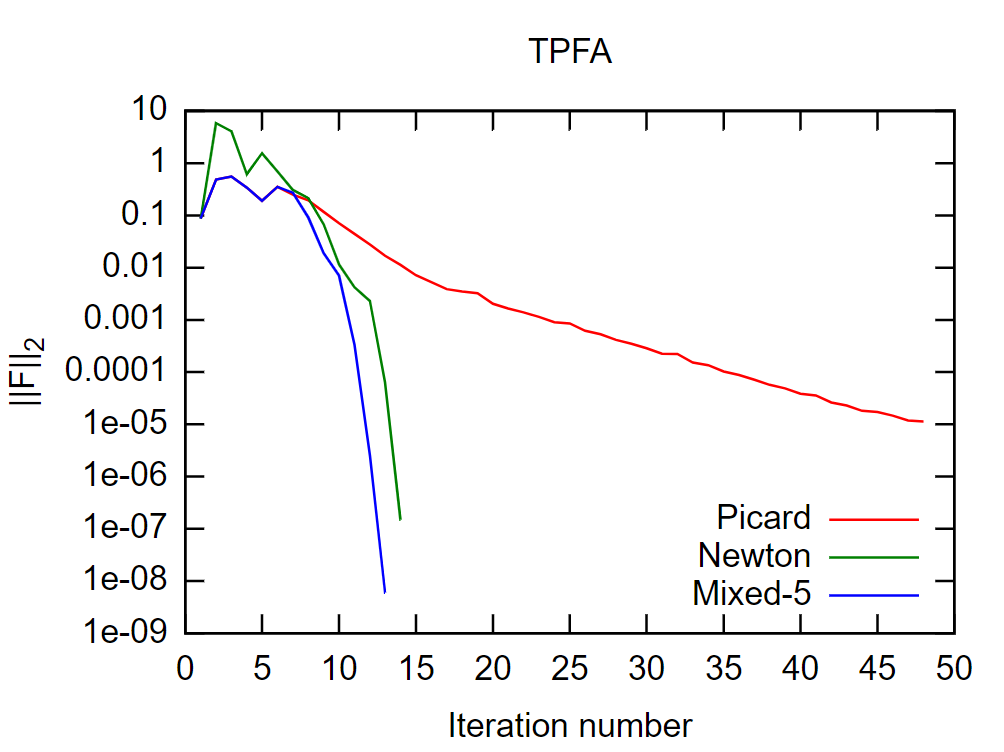}
\includegraphics[width=0.4\textwidth]{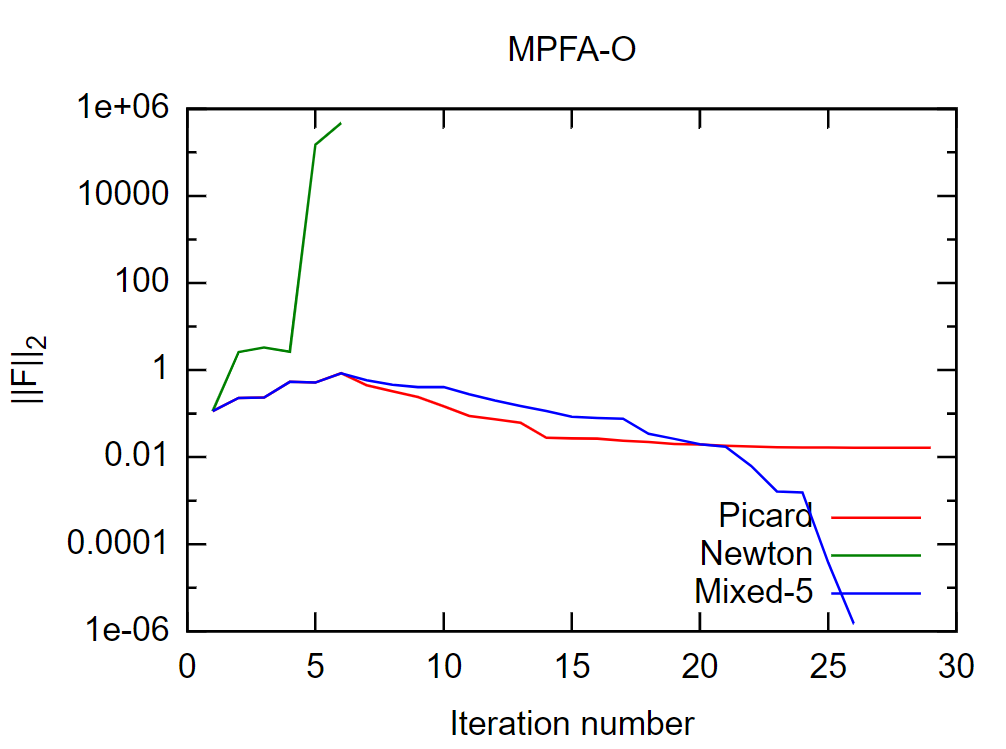}
\includegraphics[width=0.4\textwidth]{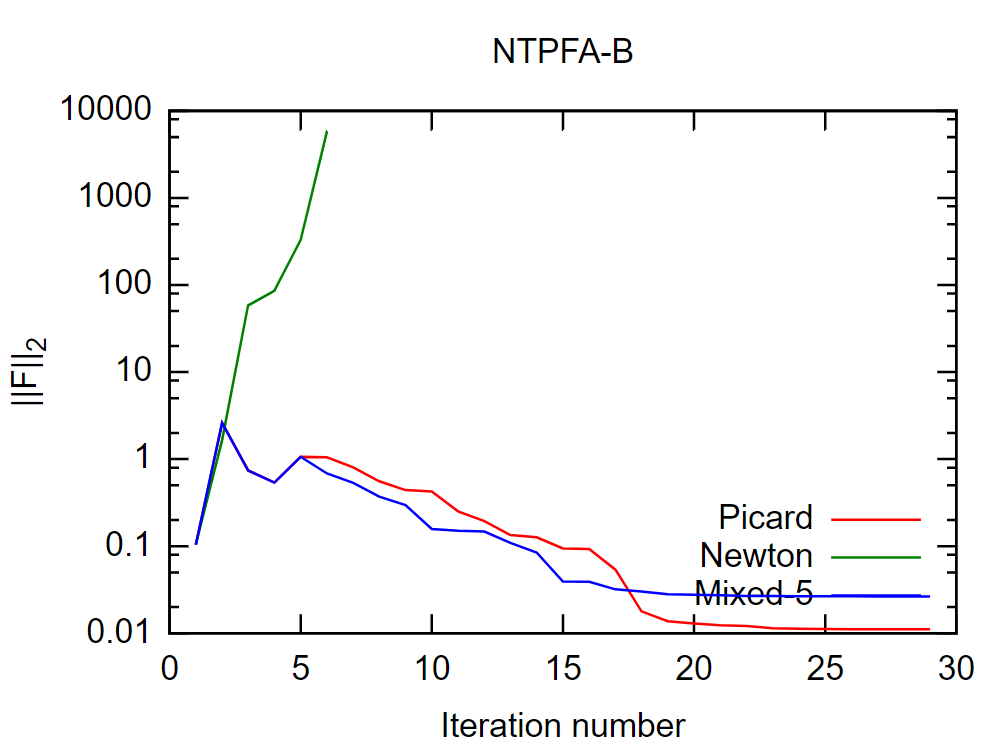}
\includegraphics[width=0.4\textwidth]{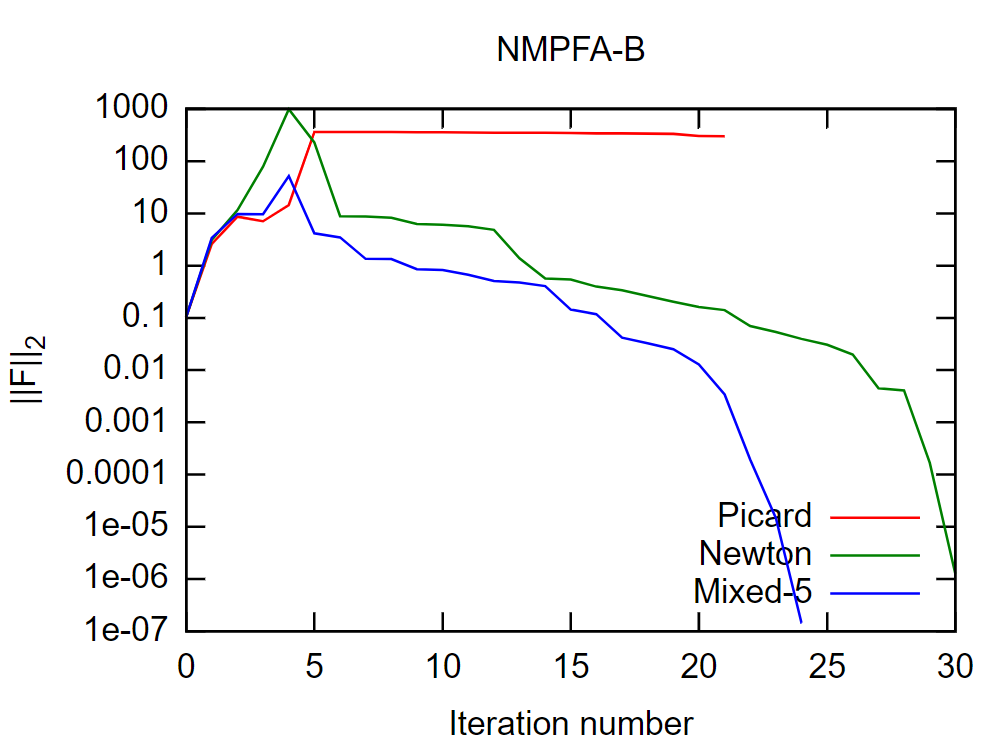}
\caption{Convergence for the modified dam problem on Cartesian grid, 6400 cells}\label{pic:dam_rot_6400}
\end{figure}

 Obtained saturation distributions for NMPFA-B and TPFA are shown in figure \ref{pic:dam_rot_sat}. Note that saturation profile differs significantly for TPFA which does not give approximation on non-$\mathbb{K}$-orthogonal meshes. Comparison of convergence of three solvers: Picard, Newton and mixed solver with 5 iterations of Picard before switching to Newton are shown in figures \ref{pic:dam_rot_400} and \ref{pic:dam_rot_6400}. As expected, the Picard method has the slowest convergence and often experiences line search failures. Of the two other solvers, mixed Picard-Newton solver turned out to be more robust. On coarser grid, it was a little outperformed by the pure Newton method in case of linear FV schemes, but was significantly faster for NMPFA-B scheme. On finer grid, it was faster for TPFA and NMPFA-B schemes and was the only method to converge for MPFA-O scheme. Surprisingly, NTPFA-B case exhibited worse convergence than NMPFA-B one, with none of the solvers being able to achieve convergence. It should be also noted that pure Newton solver sometimes tended to diverge at the beginning, and mixed Picard-Newton solver negates this drawback. 
	\par
	\subsection{Unconfined flow model: real site}
	As mentioned already, use of the unconfined flow model for real site results in relatively easy convergence: in some cases, the nonlinearity continuation method is able to converge in 1 step (from $q = 0$ directly to $q = 1$). Therefore, this scenario could be the ultimate goal for nonlinear solver. It should be mentioned that since the expressions in unconfined flow model are mesh-dependent, convergence behavior may significantly differ for various meshes, and the following tests are presented for grids leading to relatively good convergence (mostly, these grids are quite coarse in the vertical direction). Still, some solvers struggle in these simpler cases. 
	\par
\begin{figure}
		\centering
		\includegraphics[width=0.6\textwidth]{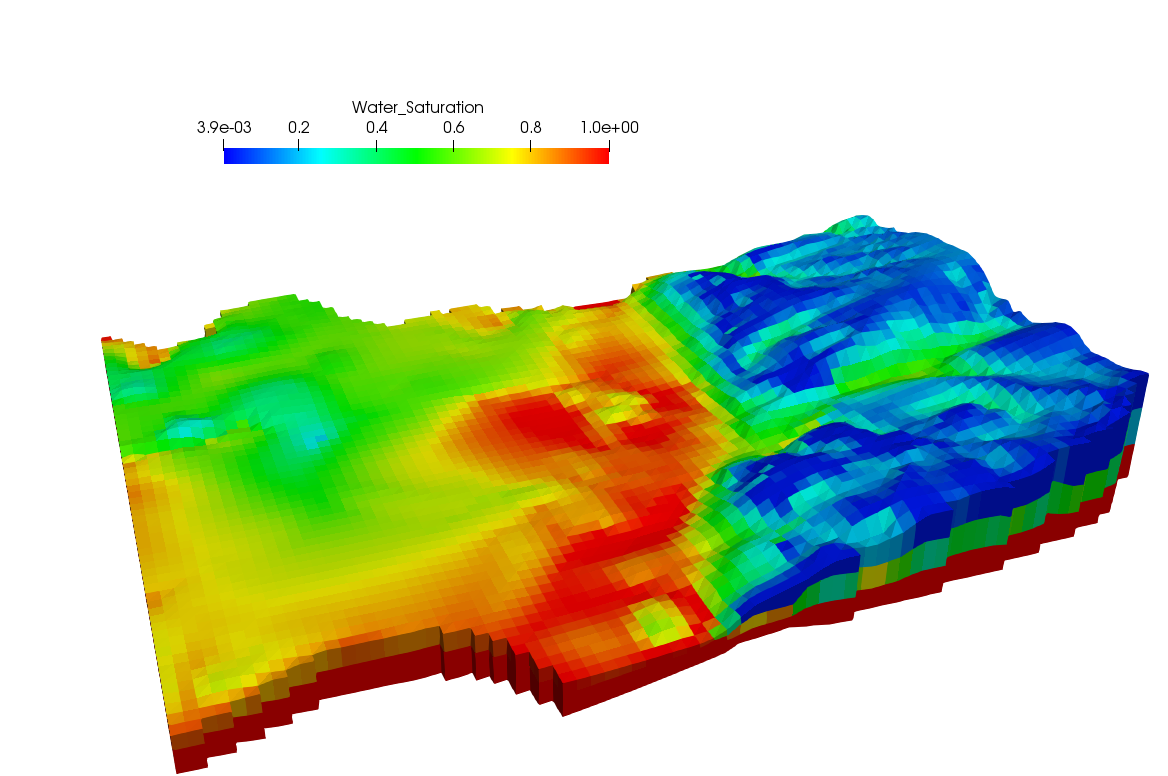}
		\includegraphics[width=0.6\textwidth]{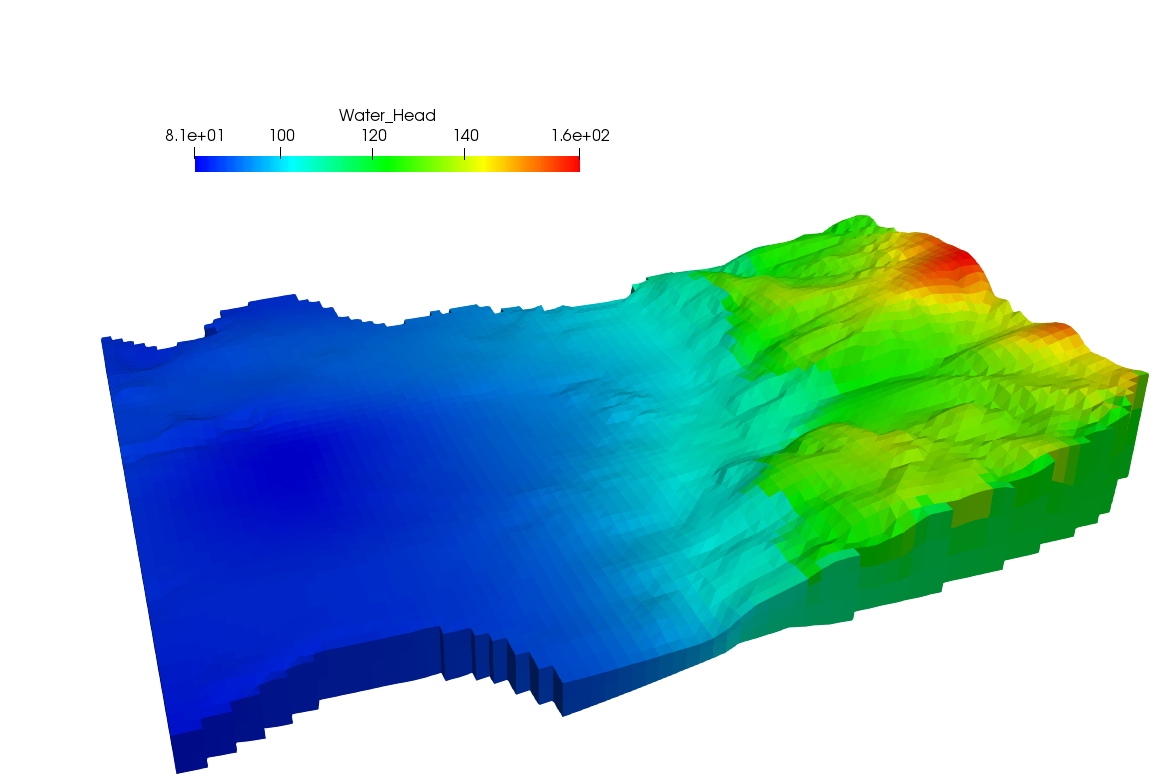}
		\caption{Saturation and hydraulic head distribution, real site, 13000 cells, NTPFA-B}\label{pic:shk_unconf_13k_pic}
	\end{figure}
\begin{figure}
		\centering
		\includegraphics[width=0.4\textwidth]{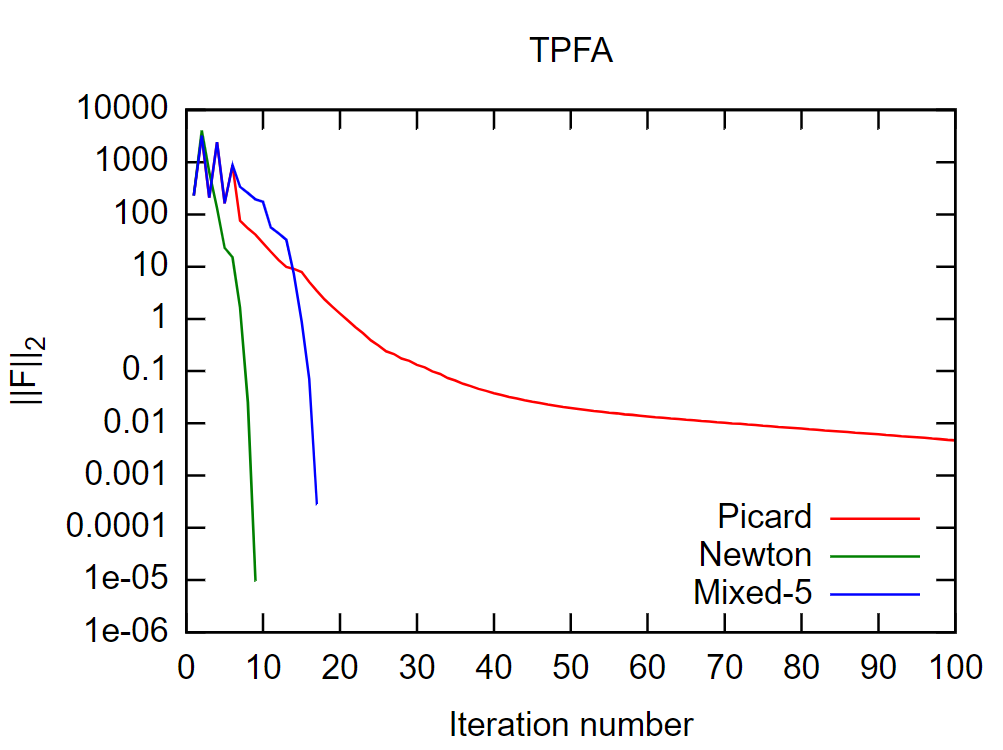}
		\includegraphics[width=0.4\textwidth]{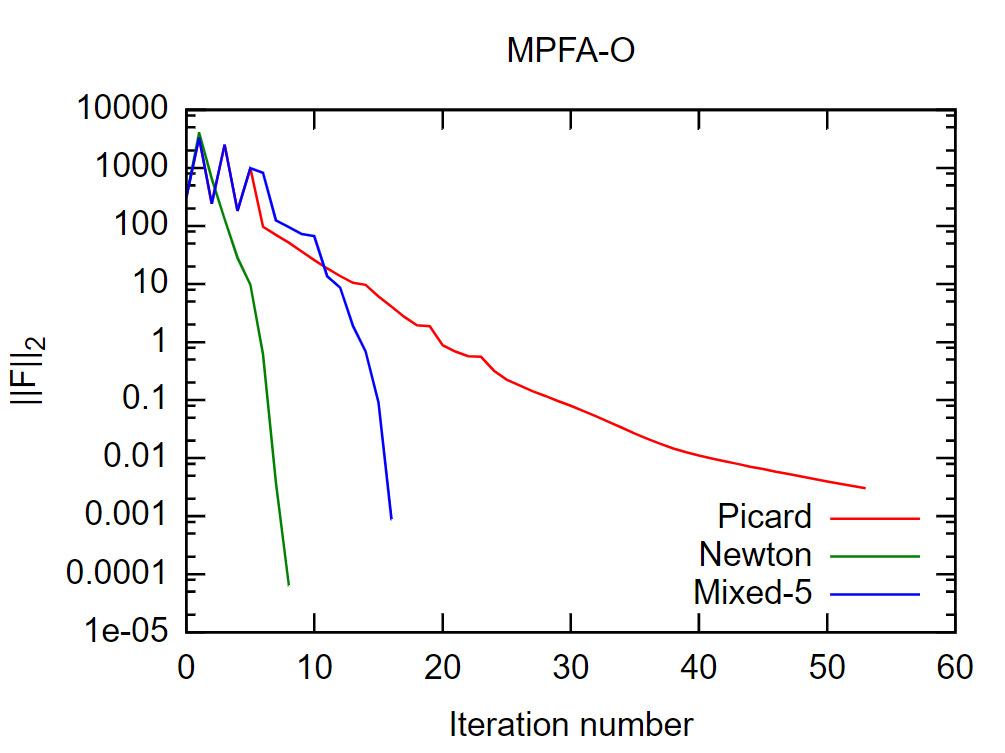}
		\includegraphics[width=0.4\textwidth]{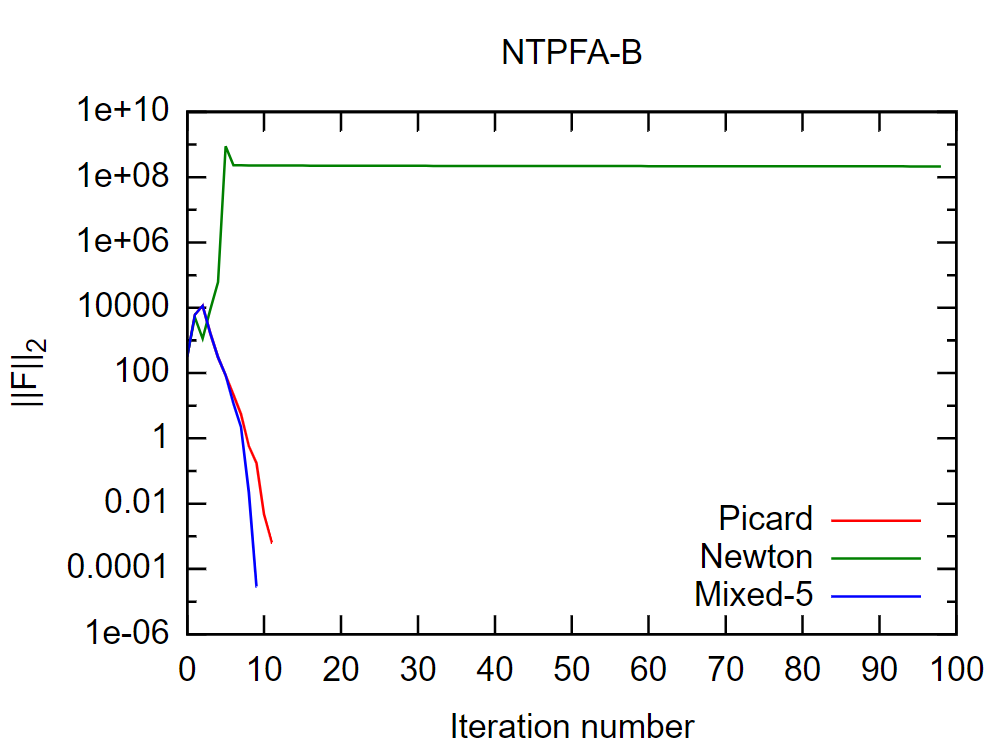}
		\includegraphics[width=0.4\textwidth]{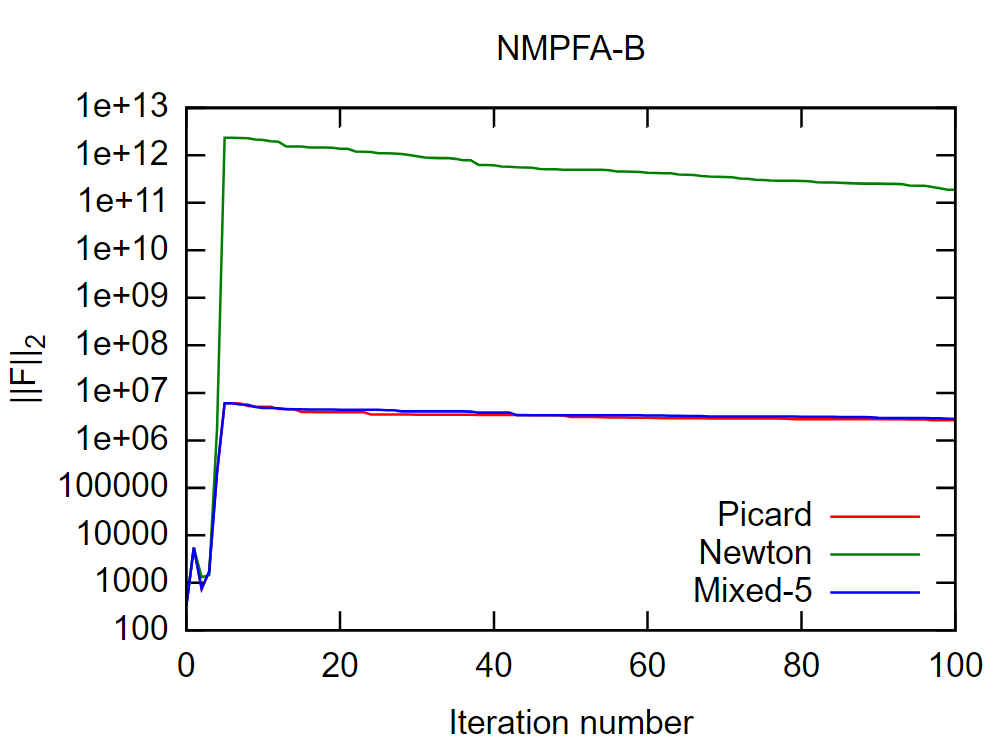}
		\caption{Convergence for the real site, unconfined flow model, 13000 cells}\label{pic:shk_unconf_13k}
	\end{figure}

	The test was conducted for a 13000-cell hexahedral mesh and the following solver parameters: $\varepsilon_{rel} = 10^{-5},~~\varepsilon_{abs} = 10^{-6},~~nit_{\max} = 100$. Obtained saturation and head distributions are shown in figure \ref{pic:shk_unconf_13k_pic}. Convergence plot are presented in figure \ref{pic:shk_unconf_13k}. Results showed that for linear FV schemes, TPFA and MPFA-O, almost all solvers worked successfully except for Picard in TPFA case. Newton solver had no convergence problems and therefore was faster than the mixed solver. For the NTPFA-B scheme, Newton method was unable to converge, mixed solver again was the fastest and, surprisingly, Picard solver exhibited fast convergence. For the NMPFA-B scheme, none of the solvers were able to converge. The reason of this failure is the sharp increase in $||F||$ in the first iterations. This led to the use of technique described in section \ref{sec:LSS}: use of relaxation with fixed small $\omega$ at several first iterations, which is likely to prevent initial sharp increase in $||F||$. This technique was used for the mixed solver, since this solver has shown superior performance already, and led to success. Only 5 initial Picard iterations with fixed $\omega = 0.1$ were enough to overcome initial sharp rise of $||F||$ and change overall convergence for NMPFA-B.
	    
\begin{figure}
\centering
\includegraphics[width=0.8\textwidth]{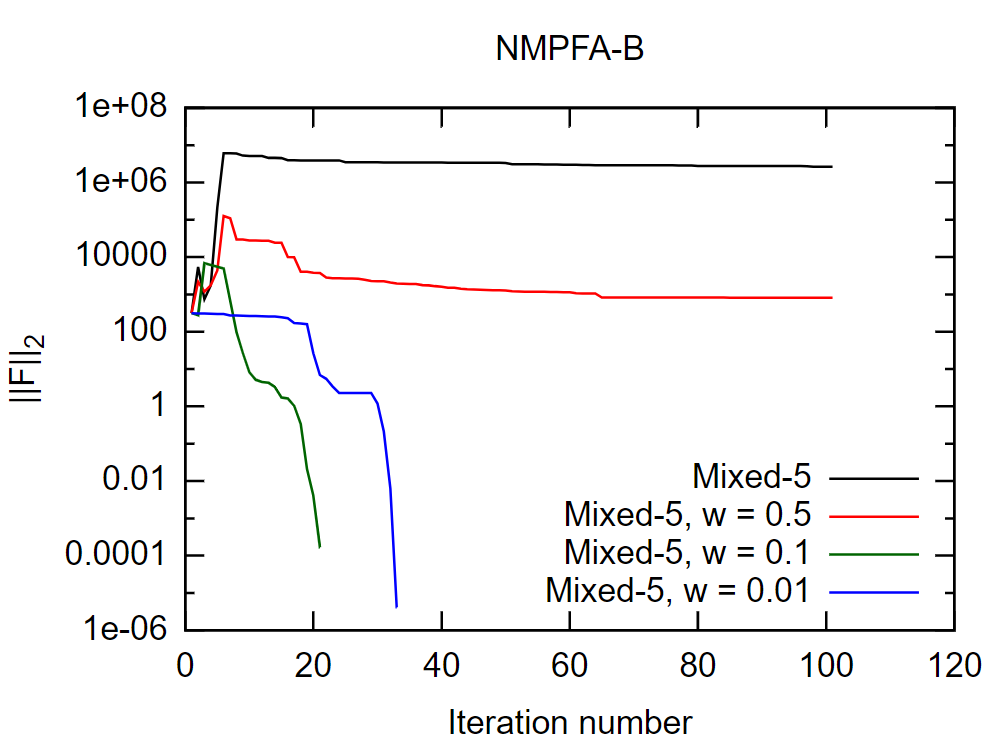}
\caption{Peformance of the mixed solver with fixed relaxation parameter $\omega$ at first 5 iterations, real site, unconfined flow model, 13000 cells}\label{pic:shk_unconf_13k_fix}
\end{figure}

    \subsection{Van Genuchten -- Mualem model: modified dam test}
    The purpose of this test is to examine solvers' performance on a simple problem in case of more nonlinear van Genuchten and Mualem expressions \eqref{eq:vgm_theta} -- \eqref{eq:vgm_Kr}. It is known that relative permeability \eqref{eq:vgm_Kr} has discontinuity for $n < 2$ \cite{lipnikov2016new}, and in this case such $n = 1.2  $ was chosen. Strong nonlinearity of the constitutive relationships, rotated anisotropic tensor and use of non-$\mathbb{K}$-orthogonal meshes still make this problem quite hard to solve, requiring several continuation steps. Therefore, it is a good test on how various solvers affect the number of continuation steps
    \par
\begin{figure}
    	\centering
    	\includegraphics[width=0.4\textwidth]{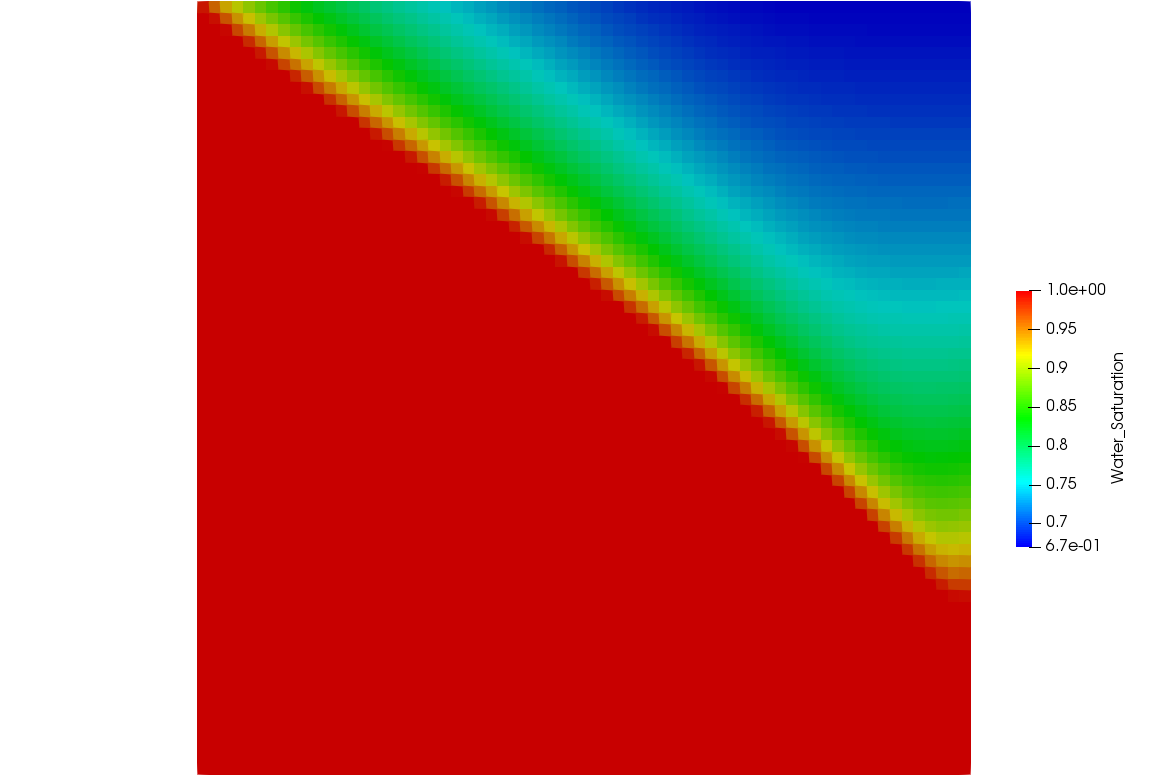}
    	\includegraphics[width=0.4\textwidth]{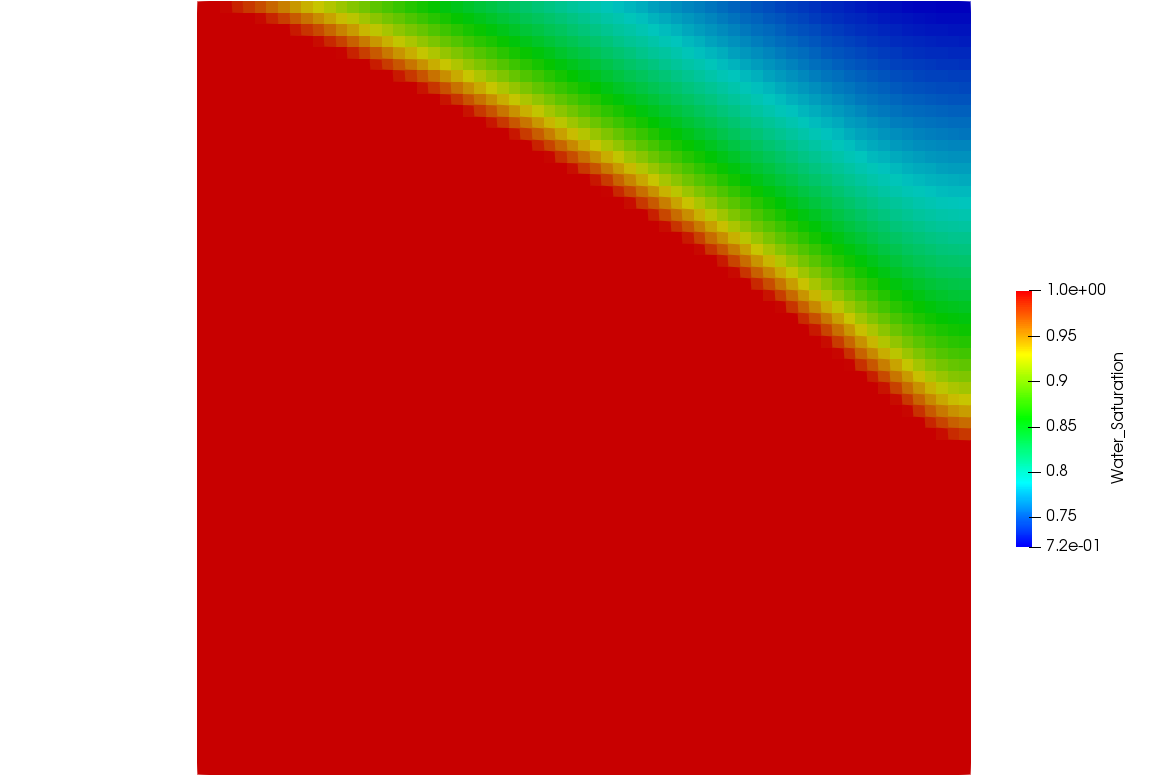}
    	\caption{Saturation distributions produced by NMPFA-B and TPFA for the modified dam problem, van Genuchten -- Mualem model, 6400 cells}\label{pic:dam_rot_sat_vgm}
    \end{figure}

    The following solver parameters: $\varepsilon_{rel} = 10^{-5},~~\varepsilon_{abs} = 10^{-5},~~nit_{\max} = 80$ were used. Tests were carried out on several grids including both triangular and Cartesian. Results for Cartesian mesh of 5500 cells are given in table \ref{tab:vgm_dam_hex} and results for triangular mesh of 1800 cells are given in table \ref{tab:vgm_dam_tri}. Obtained saturation distributions for Cartesian mesh are shown in figure \ref{pic:dam_rot_sat_vgm}. In general, the best convergence was observed for TPFA and, surprisingly, NMPFA-B schemes, in some cases it was possible to finish continuation in 1 step even with the Newton solver. NTPFA-B exhibited worse convergence compared to NMPFA-B, often requiring more than 1 continuation step when NMPFA-B required 1. Interestingly enough, for NTPFA-B on triangular mesh the mixed Picard-Newton solver reduced the number of continuation steps, but increased total iteration number and computation time. The mixed Picard-Newton solver also often lead to failure of continuation procedure with linear continuation function $\mathcal{K}_{lin}$ for TPFA and especially MPFA-O. These irregularities likely result from non-smoothness of the relative permeability function $K_r$.
    \par 
   The main interest in this test was to check if advanced solver, e.g. mixed solver with 5 Picard iterations with fixed $\omega = 0.1$ in the beginning, would reduce number of continuation steps. Tests gave no clear answer. The mixed solver could worsen performance in cases where Newton works well and use of other solvers is redundant. This happened mostly for linear FV schemes, while for nonlinear schemes the mixed solver could reduce the number of continuation steps. However, this behavior differed on various meshes, likely due to non-smoothness of $K_r$.

\begin{table}
	\centering
	\begin{tabular}{|p{4.5cm}|r|r|r|r|r|r|}
		\hline
		\multirow{2}{*}{}
		& \multicolumn{3}{c|}{Newton} & \multicolumn{3}{c|}{Mixed Picard-Newton}\\
		%\hline
		% \textbf{Inactive Modes} & \textbf{Description}\\
		%\cline{2-5}
		& Time, s & Cont.st. & Tot.iter. & Time, s & Cont.st. & Tot.iter.\\
		\hline
		TPFA,    $\mathcal{K}_{pow}= \mathcal{K}_{lin}$ &   5.6 &  1 &  22 &  13.4 &  1 &  53 \\
		MPFA-O,  $\mathcal{K}_{lin}$                    &  74.2 &  1 &  25 & Fail  & 15 &1321 \\
		MPFA-O,  $\mathcal{K}_{pow}$                    &  74.2 &  1 &  25 &3624.0 &  5 & 297 \\
		NTPFA-B, $\mathcal{K}_{lin}$                    &  52.4 &  2 &  61 &  12.9 &  1 &  29 \\
		NTPFA-B, $\mathcal{K}_{pow}$                    &  43.1 &  2 &  62 &  12.9 &  1 &  29 \\
		NMPFA-B, $\mathcal{K}_{lin}$                    &  96.3 &  2 &  90 &  12.9 &  1 &  29 \\
		NMPFA-B, $\mathcal{K}_{pow}                   $ &  79.2 &  2 & 106 &   7.8 &  1 &  40 \\\hline
	\end{tabular}
	\caption{Comparison of solvers for the modified dam test with van Genuchten -- Mualem functions, Cartesian grid with 5500 cells}
	\label{tab:vgm_dam_hex}
\end{table}

\begin{table}
	\centering
	\begin{tabular}{|p{4.5cm}|r|r|r|r|r|r|}
		\hline
		\multirow{2}{*}{}
		  & \multicolumn{3}{c|}{Newton} & \multicolumn{3}{c|}{Mixed Picard-Newton}\\
		 %\hline
		% \textbf{Inactive Modes} & \textbf{Description}\\
		%\cline{2-5}
		        & Time, s & Cont.st. & Tot.iter. & Time, s & Cont.st. & Tot.iter.\\
		\hline
		TPFA,    $\mathcal{K}_{lin}                   $ &   5.1 &  1 &  29 & Fail  & 55 &2106 \\
		TPFA,    $\mathcal{K}_{pow}                   $ &   5.1 &  1 &  29 &  52.1 &  9 & 341 \\
		MPFA-O,  $\mathcal{K}_{lin}$                    &  28.2 &  1 &  25 & Fail  & 12 & 857 \\
		MPFA-O,  $\mathcal{K}_{pow}$                    &  28.2 &  1 &  25 & 729.0 & 14 & 495 \\
		NTPFA-B, $\mathcal{K}_{lin}$                    &   4.1 &  2 &  26 &   5.4 &  1 &  34 \\
		NTPFA-B, $\mathcal{K}_{pow}$                    &  29.4 &  3 & 109 &   5.4 &  1 &  34 \\
		NMPFA-B, $\mathcal{K}_{pow}=\mathcal{K}_{lin} $ &   2.2 &  1 &  14 &   7.8 &  1 &  40 \\\hline
	\end{tabular}
	\caption{Comparison of solvers for the modified dam test with van Genuchten -- Mualem functions, triangular mesh with 1900 cells}
	\label{tab:vgm_dam_tri}
\end{table}
	
	\section{Conclusion}
	Different iterative solvers: Newton, Picard and mixed Picard-Newton solver, all combined with line search technique based on Armijo rule, were tested within the nonlinearity continuation method for the steady-state Richards equation. Slow convergence of the Picard method makes it unattractive as a standalone solver, but it is useful in the mixed solver as a way to improve initial guess. Another way to improve initial convergence is to apply relaxation with fixed parameter at first iterations. 
	\par The solvers were tested on systems arising in finite volume discretization on unstructured grids. Finite volume schemes with linear and nonlinear, two- and multi-point flux approximations were used. Tests showed that for problems with mild nonlinearity it is often possible to finish nonlinearity continuation in only 1 step by use of the mixed solver with fixed relaxation parameter at first iterations. For problems with stronger nonlinearity, the mixed solver behaved differently: it may be helpful for schemes with nonlinear flux approximation, but often deteriorated performance for the schemes with linear flux approximation. In some cases, the mixed solver decreased number of continuation steps while increasing the total iteration number.
	
	\section*{Acknowledgements}
	The reported study was funded by Russian Foundation for Basic Research (RFBR), project number 20-31-90126.
	
	%% If you have bibdatabase file and want bibtex to generate the
	%% bibitems, please use
	%%
	\bibliographystyle{elsarticle-num} 
	\bibliography{ljm2021.bib}
	
	%% else use the following coding to input the bibitems directly in the
	%% TeX file.
	
	%\begin{thebibliography}{00}
	
	%% \bibitem[Author(year)]{label}
	%% Text of bibliographic item
	
	%\bibitem{ReddyGartling}
	%Reddy J. N., Gartling D. K. (2010) The finite element method in heat transfer and fluid %dynamics. – CRC Press.
	
	%\end{thebibliography}
\end{document}